\documentclass[12pt,a4paper,draft]{article}
\usepackage[english,russian]{babel}
\usepackage{amssymb, amsmath, latexsym, amsfonts, amsthm}
\begin{document}
\renewcommand{\refname}{References}
\newtheorem*{theorem}{Theorem}
\newtheorem*{theorem M}{Theorem M}
\newtheorem*{lemma}{Lemma}
\newtheorem{corollary}{Corollary}
\newtheorem{proposition}{Proposition}
\begin{center}
On perturbation of a surjective convolution operator
\end{center}
\begin{center}
Il'dar Kh. Musin\footnote {E-mail: musin\_ildar@mail.ru}
\end{center}
\begin{center}
Institute of Mathematics with Computer Centre of Ufa Scientific Centre of Russian Academy of Sciences, 

Chernyshevsky str., 112, Ufa, 450077, Russia
\end{center}

\vspace {0.3cm}

\renewcommand{\abstractname}{}
\begin{abstract}
{\bf Abstract:} 
Let $\mu \in {\cal E}'({\mathbb R}^n)$ be a compactly supported distribution such that its support is a convex set with non-empty interior.
Let $X_2$ be a convex domain in ${\mathbb R}^n$, $X_1 = X_2 + supp \   \mu $. 
Assuming that a convolution  operator 
$A: {\cal E}(X_1) \to {\cal E}(X_2)$ acting by the rule
$(Af)(x) = (\mu * f)(x)
$
is surjective we provide a condition on a linear continuous operator $B: {\cal E}(X_1) \to {\cal E}(X_2)$ that guarantees surjectivity of the operator $A+B$.

\vspace {0.3cm}


\vspace {0.3cm}


\end{abstract}

\vspace {1cm}

\begin{section}
{\bf  Introduction} 
\end{section}

{\bf 1.1. Statement of the problem and the main result}.  
By ${\cal E}(X)$, where $X$ is an open subset of ${\mathbb R}^n$, we mean 
the space of infinitely differentiable functions in  $X$ with the topology defined by the system of semi-norms 
$$
{\Vert f \Vert}_{K, N} =
\displaystyle \sup \limits_{x \in K, \vert \alpha \vert \le N}
\vert (D^{\alpha}  f)(x) \vert, \ K \Subset X, N \in {\mathbb Z}_+.
$$
Its strong dual space ${\cal E}'(X)$ is the space of distributions with compact support in $X$. 

If $0 \ne \mu \in {\cal E}'({\mathbb R}^n)$ and $X_1$, $X_2$ are two non-empty open sets in ${\mathbb R}^n$ such that
\begin {equation}
\label{s0}
X_2 + supp \ \mu \subset X_1,
\end {equation}
then the convolution $\mu * f$ of distribution $\mu$ and a function $f \in {\cal E}(X_1)$ defined by the rule 
$$
(\mu * f)(x) = \mu(f(x +y)), \ x \in X_2,
$$
is in ${\cal E}(X_2)$.

The problem of characterizing the surjectivity of convolution operators and, in particular, of partial differential operators, has interested several authors. L. Ehrenpreis \cite {E54} and B. Malgrange \cite {Malgr} established that for each nonzero polynomial $P$ of $n$ variables $P(D)(C^{\infty}({\mathbb R}^n)) = C^{\infty}({\mathbb R}^n)$.
For a distribution $\mu \in {\cal E}'({\mathbb R}^n),\mu \ne 0$, 
L. Ehrenpreis \cite {E} proved that the convolution operator $f \to \mu * f$ acting from 
$C^{\infty}({\mathbb R}^n)$ to $C^{\infty}({\mathbb R}^n)$ is surjective if and only if $\mu$ is invertible, respectively,  if
its Fourier-Laplace transform $\hat \mu$ defined by 
$$
\hat \mu (z) = \mu(e^{\langle -iz, \xi \rangle}), \  z \in {\mathbb C}^n,
$$
is slowly decreasing, that is, if 
there exists a constant $a > 0$ such that for every 
$\xi \in {\mathbb R}^n$ there exists $\eta \in {\mathbb R}^n$  satisfying 
$$
\Vert \xi - \eta \Vert 
\le a \ln (2 + \Vert \xi \Vert) 
$$
and
$$
\vert \hat \mu (\eta) \vert
\ge (a + \Vert \xi \Vert)^{-a}.
$$
The answer to a question concerning the surjectivity of convolution operators in general case
was given by L. H\"ormander \cite {H}-\cite {H2}. 
He proved that the convolution equation 
$
\mu * f = g
$
has a solution $f \in {\cal E}(X_1)$ for every $g \in {\cal E}(X_2)$ 
if and only if $\mu$ 
is invertible and the pair $(X_1, X_2)$  is $\mu$-convex for supports.
Recall that the pair $(X_1, X_2)$ of open sets
$X_1, X_2$ in ${\mathbb R}^n$ satisfying (\ref {s0}) 
is called $\mu$-convex for supports \cite [Definition 3.2] {H}, \cite [Definition 3.2] {H1}
if for any $\nu \in {\cal E}'(X_2)$ we have
$$
dist(supp \ \nu,
{\mathbb R}^n \setminus X_2) =
dist(supp \ \mu * \nu, {\mathbb R}^n \setminus X_1).
$$
Here $\mu*\nu$ is the convolution of distributions $\mu$ and $\nu$ defined by 
$$
(\mu*\nu) (f) = \mu(\nu(f(x + y))), \ f \in {\cal E}(X_1),
$$
$dist (A, B) = \inf \{\Vert x - y \Vert: x \in A, y \in B\}, \Vert \cdot \Vert$ 
is the Euclidean norm in ${\mathbb R}^n$.

It was proved by L. H\"ormander \cite [Theorem 5.4, Corollary 5.4] {H3} that if $\mu_1, \mu_2 \in {\cal E}'({\mathbb R}^n)$ have disjoint singular supports and $\mu_1$ is slowly decreasing then $\mu_1 + \mu_2$  is also slowly decreasing. Later a direct proof of this H\"ormander's result was given by W. Abramczuk \cite [Theorem 1]{A}. 
Thus, if $\mu_1, \mu_2 \in {\cal E}'({\mathbb R}^n)$ have disjoint singular supports and $\mu_1$ defines a surjective convolution operator on ${\cal E}({\mathbb R}^n)$ then the convolution operator associated to $\mu_1 + \mu_2$ is also surjective. 

Afterwards there were not many works devoted to perturbations of convolution operators in spaces of infinitely differentiable functions.  Among them it is necessary to mention relatively recent work of C. Fernandez, A. Galbis and D. Jornet \cite {FGJ} where the behavior of the perturbed convolution operator on spaces of ultradifferentiable functions in the sense of Braun, Meise and Taylor \cite {BMT} was studied. They essentially used results of  J. Bonet, A. Galbis and R. Meise \cite {BGM} on the range of convolution operators on spaces of non-quasianalytic ultradifferentiable functions and results of R. Braun, R. Meise, D. Vogt \cite {BMV} on surjectivity of convolution
operators on classes of ultradifferentiable functions. 

In the present note the problem of surjectivity of perturbed convolution operators is studied in spaces of infinitely differentiable functions on convex domains of ${\mathbb R}^n$. 
The statement of the problem differ from ones in \cite {H3}, \cite {A}. 
It is inspired by researches of S.G. Merzlyakov \cite {Mer} of perturbations of convolution operators in spaces of holomorphic functions.
Namely, fix a distribution $\mu \in {\cal E}'({\mathbb R}^n)$ 
such that its support is a convex set with non-empty interior. 
Let $X_2$ be a convex domain in ${\mathbb R}^n$, $X_1 = X_2 + supp \ \mu $. 
Note that in this case the pair $(X_1, X_2)$ is $\mu$-convex for supports. 
It follows from theorem on supports \cite[Theorem 4.3.3]{H2}
and from the fact that for arbitrary convex domain 
$\Omega \subset {\mathbb R}^n$ and for each compact 
$K \subset \Omega$ we have that
$
dist(K, {\mathbb R}^n \setminus \Omega) = dist(ch K, {\mathbb R}^n \setminus \Omega).
$
Here $ch K$ is a convex envelope of a compact $K$.
Assume that the convolution  operator 
$A: {\cal E}(X_1) \to {\cal E}(X_2)$ acting by the rule
$(Af)(x) = (\mu * f)(x)
$
is surjective (thus, $\mu$ is invertible). 
Consider a linear operator $B: {\cal E}(X_1) \to {\cal E}(X_2)$ such that for any convex compact 
$K_2$ in $X_2$ there exist a convex compact subset $V$ of the interior of support of $\mu$ (denoted by $supp \ \mu$) and a number 
$N_1 \in {\mathbb Z}_+$ such that for each positive $\varepsilon$ that is less than the distance between $K_2$ and the boundary of $X_2$ 
and for each $N_2 \in {\mathbb Z}_+$ there exists a number
$c = c(\varepsilon, N_2) > 0 $ such that
\begin {equation}
\label{ss1}
{\Vert Bf \Vert}_{K_2^{\varepsilon}, N_2} \le
c {\Vert f \Vert}_{K_2 + V, N_1}, \
f \in {\cal E}(X_1).
\end {equation}

The main result of the paper is the following

\begin{theorem}
The operator $A+B: {\cal E}(X_1) \to {\cal E}(X_2)$ is surjective.
\end{theorem}

{\bf 1.2. Organization of the article}. In  section 2 there are given two useful auxiliary results. 
The first is the Phragmen-Lindel\"of type result (see Proposition 1). The second is the Division Theorem of L. H\"ormander (see e.g., \cite [Corollary 2.6.] {H4}). Also here we recall definitions of two types of locally convex spaces introduced by Jos\'e Sebasti\~ao e Silva in \cite {S-S}. The main result is proved in section 3. In section 4 we give an example of the operator $B$. 

{\bf 1.3. Some notations}. For $u=(u_1, \ldots , u_n) \in {\mathbb R}^n \ ({\mathbb C}^n)$, $v=(v_1, \ldots , v_n) \in {\mathbb R}^n \ ({\mathbb C}^n)$ \ 
$\langle u, v \rangle  = u_1 v_1 + \cdots + u_n v_n$ and $\Vert u \Vert$ denotes the Euclidean norm in ${\mathbb R}^n ({\mathbb C}^n)$. 

For $\alpha = (\alpha_1, \ldots , \alpha_n) \in {\mathbb Z}^n$ \ \  $\vert \alpha \vert = \alpha_1 + \ldots  + \alpha_n$, $D^{\alpha}$ is the corresponding derivative.

If $\Omega \subset {\mathbb R}^n$ then $\overline \Omega$, 
$int \ \Omega$, $\partial \Omega$, $ch \ \Omega$ denote its closure, interior, boundary and convex envelope, respectively. For $\varepsilon > 0$ let $\Omega^{\varepsilon} = \{x \in {\mathbb R}^n: 
\Vert x - y \Vert \le \varepsilon \ \text {for some} \ y \in \Omega \}$.

For $r>0$ let $D(r) =\{x \in {\mathbb R}^n: \Vert x \Vert < r \}$.

Supporting function $H_K$ of convex compact set $K \subset {\mathbb R}^n$ is defined by
$
H_K(y)=\displaystyle \max_{t \in K} \langle y, t \rangle, \ y \in {\mathbb R}^n.
$

$H({\mathbb C}^n)$ is the space of entire functions on ${\mathbb C}^n$. 

For a locally convex space $E$ let $E^*$ be the strong dual space.

\begin{section}
{\bf  Preliminaries} 
\end{section}

{\bf 2.1. Auxiliary results}. In the proof of the Theorem the following two results will be useful.

\begin{proposition}
\label{FL}
Let $b$ be a non-negative convex positively homogeneous of order 1 function on ${\mathbb C}^n$ and $g \in H({\mathbb C}^n)$. Assume that for each $\varepsilon > 0$ there exists a constant
$c_{\varepsilon} > 0$ such that
$$
\vert g(z) \vert \le
c_{\varepsilon}
\exp {(b(z) + {\varepsilon} \Vert z \Vert)}, \
z \in {\mathbb C}^n,
$$
and for some $M >0$ and $N \in {\mathbb Z}_+$
$$
\vert g(x) \vert \le M (1 + \Vert x \Vert)^N, \ x \in {\mathbb R}^n.
$$

Then
$$
\vert g(z) \vert \le
2^{\frac {N}{2}} M (1 + \Vert z \Vert)^{2N} \exp (b(iIm \ z)), \
z \in {\mathbb C}^n.
$$
\end{proposition}

It is an easy consequence of the Lemma below that was proved in fact in 
\cite {MusSpain}. To formulate it let us introduce a space ${\cal P}_a(T_C)$ as follows. Let $C$ be an open convex cone in ${\mathbb R}^n$ with an apex at the origin
and $a$ is a nonnegative convex contionuous positively homogeneous function of degree 1 on ${\mathbb R}^n + i {\overline C}$. Then ${\cal P}_a(T_C)$ is the space of functions $f$ holomorphic on tube domain $T_C = {\mathbb R}^n + i C$ and satisfying the condition: for any $\varepsilon > 0$ there exists a constant $c = c_{\varepsilon, f} > 0$ such that
$$
\vert f(z) \vert \le
c \exp {(a(z) + {\varepsilon} \Vert z \Vert)}, \ z \in {\mathbb R}^n + i C.
$$

\begin{lemma}
Let $g \in {\cal P}_a(T_C)$ and for $\xi \in {\mathbb R}^n$ \
$
{\varlimsup  \limits_{z \to \xi, \atop z \in T_C}} \vert g(z) \vert \le M.
$

Then 
$$
\vert g(x + i y) \vert \le M \exp (a(iy)), \
x + i y \in T_C.
$$
\end{lemma}

{\bf Remark}. In \cite [Lemma] {MusSpain} it is assumed that $C$ is acute. Analysis of the proof of this Lemma shows that this condition on $C$ is unnecessary. 

The following result was obtained by L. H\"ormander (see e.g., \cite [Corollary 2.6.] {H4}).
\begin{proposition}
For j= 1, 2, 3 let $u_j \in {\cal E}'({\mathbb R}^n)$, let
$$
H_j(\eta) = sup \{\langle x, \eta \rangle, \ x \in  supp \ u_j\},
$$
and let $U_j$ be the Fourier-Laplace transform of $u_j$. Assume that $U_2 = \frac {U_3} {U_1}$ 
is entire. Then it follows that $H_2 = H_3-H_1$ is a supporting function and
that for every $\varepsilon > 0$
$$
\vert U_2 (\zeta) \vert \le 
C_{\varepsilon} \exp (H_2 (Im \zeta) + \varepsilon \Vert \zeta \Vert), \ 
\zeta \in {\mathbb C}^n.
$$

\end{proposition}

{\bf 2.2. Two definitions}. 
Recall the definitions of $(M^*)$-space and $(LN^*)$-space from \cite {S-S}.

{\bf Definition 1.} $(M^*)$-space is a locally convex space $F$ which is the projective limit of a sequence of normed  spaces $F_k$ with linear continuous mappings $g_{mk}: F_k \to F_m$,  $m < k,$ such that $g_{k, k+1}$  is compact for each $k \in {\mathbb N}$.

{\bf Definition 2.} $(LN^*)$-space is a locally convex space $E$ which is the inductive limit of an increasing sequence of normed spaces $E_k$ such that the unit ball of $E_k$ is relatively compact in $E_{k+1}$ for each $k \in {\mathbb N}$, i.e. such that the inclusion map from $E_k$ into $E_{k+1}$ is compact.

{\bf 2.3. Some additional notations and notions used in the proof of Theorem}. 
If $X$ is an open set in ${\mathbb R}^n$ and $(K_m)_{m=1}^{\infty}$ is a sequence of compact subsets of $X$ such that 
$K_m \subset int \ K_{m+1}$ ($m =1, 2, \ldots )$ and $X = \cup_{k=1}^{\infty} K_m$ 
then let $C^m(K_m)$ be a normed space of functions $f$ smooth up to the order $m$ in $K_m$
with a norm 
$$
p_m(f) =
\sup_{x \in K_m, \vert \alpha \vert \le m}
\vert (D^{\alpha} f)(x) \vert .
$$
Note that  ${\cal E}(X)$ is a projective limit of spaces $C^m(K_m)$. 
Moreover, ${\cal E}(X)$ is dense in each  $C^m(K_m)$ and embeddings $i_m: C^{m+1}(K_{m+1}) \to C^m(K_m)$ are compact. So ${\cal E}(X)$ is an $(M^*)$-space.
Hence, ${\cal E}^*(X)$ is $(LN^*)$ space and ${\cal E}^*(\Omega)$ is an inductive limit of spaces $(C^m(K_m))^*$ \cite [Theorem 5] {S-S}.

\begin{section}
{\bf  Proof of the Theorem} 
\end{section}

The theorem will be proved if we show that the image of the operator $A+B$ is closed and dense in ${\cal E}(X_2)$.

First show that the image of $A+B$ is closed in ${\cal E}(X_2)$.
Since ${\cal E}(X_1)$ and ${\cal E}(X_2)$ are Frechet spaces then closedness of the  image of the operator $A+B$ is equivalent to closedness of the image of an adjoint operator
$(A+B)^*$ \cite[8.6.13, Theorem]{Ed}. 
Since ${\cal E}^*(X_1)$ is an $(LN^*)$-space 
then to show that the image of the operator $(A+B)^*$ is closed it is sufficient to prove that the image of the operator $(A+B)^*$ is sequentially closed (see \cite [Proposition 8]   {S-S}).
So let functionals $S_k \in {\cal E}^*(X_2)$ be such that the sequence 
$((A+B)^*S_k)_{k=1}^{\infty}$ converges to $F \in {\cal E}^*(X_1)$ in ${\cal E}^*(X_1)$.

For each  $m \in {\mathbb N}$ let $X_{2, m}$ be open bounded convex subset of $X_2$ such that $\overline {X}_{2, m} \subset X_{2, m+1}$,
$X_2 = \displaystyle \cup_{m=1}^{\infty} \overline {X}_{2, m}$. 
Put $X_{1, m} = X_{2, m} + supp \ \mu$. Then $\overline {X}_{1, m} \subset X_{1, m+1}$,
$X_1 = \displaystyle \cup_{m=1}^{\infty} \overline {X}_{1, m}$.

Since ${\cal E}^*(X_1)$ is an $(LN^*)$-space  then by properties of $(LN^*)$-spaces \cite [Theorem 2, Corollary 1] {S-S} there is $p \in {\mathbb N}$ such that
functionals $F_k: = (A+B)^*S_k$ and $F$ belong to 
$(C^p(\overline {X}_{1, p}))^*$ and the sequence 
$(F_k)_{k=1}^{\infty}$ converges to $F$ in $(C^p(\overline {X}_{1, p}))^*$.
Thus, supports of functionals $F_k$ and $F$ are in $\overline {X}_{1, p}$ 
and the order of distributions $F_k$ and $F$ is not more than $p$. 

Let $2r_p:=dist (\overline {X}_{2, p}, \partial {X}_{2, p + 1})$, 
$\tilde {X_2}:={X}_{2, p} + D(r_p)$ and 
$\tilde {X_1}:= \tilde {X_2} + supp \ \mu$.
Note that $\tilde {X_1}$ and $\tilde {X_2}$  are bounded open convex sets in ${\mathbb R}^n$ and the pair $(\tilde {X_1}, \tilde {X_2})$ is $\mu$-convex for supports. 

Denote by $\tilde {A}$ a convolution  operator $f \in {\cal E}(\tilde {X_1}) \to \mu * f$. Obviously $\tilde {A}$ is acting from 
${\cal E}(\tilde {X_1})$ to ${\cal E}(\tilde {X_2})$ linearly and continuously 
and if $f \in {\cal E}(X_1)$ then $\tilde {A}f = \tilde {A}f$.
By the (earlier cited) result of L. H\"ormander \cite {H1}, \cite [Theorem 16.5.7] {H2} we have that $\tilde {A}({\cal E}(\tilde {X_1})) = {\cal E}(\tilde {X_2})$.

Next, using the inequality (2) the operator $B$ can be extended (uniquely) to a linear continuous operator $\tilde {B}$ acting from ${\cal E}(\tilde {X_1})$ to ${\cal E}(\tilde {X_2})$. Moreover, for each convex compact ${\tilde K_2} \subset \tilde {X_2}$ there exists a compact $V \subset int (supp \ \mu)$ and a number
$N_1 \in {\mathbb Z}_+$ such that for each $\varepsilon \in (0, dist ({\tilde K_2}, \partial \tilde {X_2}))$
and for each $N_2 \in {\mathbb Z}_+$ there exists a number
$c = c(\varepsilon, N_2) > 0 $
such that
$$
{\Vert \tilde {B} f \Vert}_{{\tilde K_2}^{\varepsilon}, N_2} \le
c {\Vert f \Vert}_{{\tilde K_2} + V, N_1}, \
f \in {\cal E}(\tilde {X_1}).
$$
Putting here ${\tilde K_2} = \overline {{X}_{2, p}}$
we see that
$\tilde {B}$ is a compact operator from 
${\cal E}(\tilde {X_1})$ to ${\cal E}(\tilde {X_2})$.
By Theorem 9.6.7 in \cite {Ed} the image of the operator 
$\tilde {A} + \tilde {B}$ is closed in ${\cal E}(\tilde {X_2})$.
Hence, the image of the operator $(\tilde {A} + \tilde {B})^*$ is closed in ${\cal E}^*(\tilde {X_1})$.

For each  $j \in {\mathbb N}$ let $X_{2, j} = { X_{2, p}} + D(\frac {j}{j+1} r_p)$.
Then $\tilde X_2 = \displaystyle \cup_{j=1}^{\infty} \overline {\tilde X}_{2, j}$,
$\tilde {X_1} = \displaystyle \cup_{j=1}^{\infty}
(\overline {\tilde X}_{2, j} + supp \ \mu)$.
Note that for some $m \in {\mathbb N}$ supports of functionals $F$,
$F_k = (A+B)^*S_k$
$(k=1, 2, \ldots )$ are in
$\overline {\tilde X}_{2, m} + supp \ \mu$.

Now 
take an arbitrary functional $S_k$ and show that convex envelope $W_k$ of its support is contained in $\overline {\tilde X}_{2, m+2}$.
Assume the contrary. Then there is a point $\xi \in W_k$ which is not in
$\overline {\tilde X}_{2, m+2}$.
Next, there exists a hyperplane in ${\mathbb R}^n$ dividing $\overline {\tilde X}_{2, m+2}$ and $\xi$. So we can find a point $y_0 \in {\mathbb R}^n$ such that
\begin{equation}
\label{stochka}
H_{W_k}(y_0) >
H_{\overline {\tilde X}_{2, m+2}}(y_0).
\end{equation}
Denote the order of distribution $S_k$ by $N_{2, k}$. Take $\delta_1 > 0$ so small that 
$W_k^{\delta_1} \Subset X_2$. 
Then there is a constant $a_{\delta_1, k} > 0$ such that 
$$
\vert (B^*S_k) (f)\vert = \vert S_k(Bf) \vert \le 
a_{\delta, k}{\Vert Bf \Vert}_{W_k^{\delta_1}, N_{2, k}},  \ f \in {\cal E}(X_1).
$$
By the condition on $B$ (see the inequality (\ref {ss1})) there are a convex compact 
$V \subset int (supp \ \mu)$ and a number $N_1 \in {\mathbb Z}_+$ such that for a taken small
$\delta_1 > 0$ there exists a constant $c_{\delta_1, k} > 0$ such that
$$
\vert (B^*S_k) (f)\vert \le
c_{\delta_1, k}{\Vert f \Vert}_{W_k + V, N_1}, \ f \in {\cal E}(X_1).
$$
From this we have that for all $z \in {\mathbb C}^n$ 
\begin{equation}
\label{s3}
\vert \widehat {(B^*S_k)}(z) \vert \le
c_{\delta_1, k}(1 + \Vert z \Vert)^{N_1} \exp(H_{W_k}(Im \ z) + H_V(Im \ z)).
\end{equation}
Taking into account that for some $d > 0$
$$
\label{s5}
H_V(x) \le H_{supp \ \mu}(x) - d \Vert Im \  x \Vert, \ x \in {\mathbb R}^n.
$$
we get from (4) that for all $z \in {\mathbb C}^n$ 
\begin{equation}
\vert \widehat {(B^*S_k)}(z) \vert \le
c_{\delta_1, k}(1 + \Vert z \Vert)^{N_1} 
e^{H_{W_k}(Im \ z) + H_{supp \ \mu}(Im \ z) - d \Vert Im \  z \Vert}.
\end{equation}
Further, since $F_k \in {\cal E}^*(X_1)$, 
$supp \ F_k \subset \overline {\tilde X}_{2, m} + supp \ \mu$ and  the order of distribution
$F_k$ is not more than $p$, then for each $\delta > 0$ 
there exists a constant $m_{\delta, k}>0$ such that
for each $z \in {\mathbb C}^n$
\begin{equation}
\vert \hat {F_k}(z) \vert
\le m_{\delta, k} (1 + \Vert z \vert)^p
\exp(H_{\overline {\tilde X}_{2, m}} (Im \ z) + H_{supp \ \mu}(Im \ z) +
\delta \Vert Im \  z \Vert).
\end{equation}
Using estimates (5) and (6) with $\delta = \frac {r_p}{2(m+1)(m+2)}$
we obtain that for all $z \in {\mathbb C}^n$
\begin{equation}
\vert \widehat {(A^*S_k)}(z) \vert \le a (1 + \Vert z \vert)^b
e^{
H_{ch (W_k \cup \overline {\tilde X}_{2, m+1}) + supp \ \mu} (Im \ z)
- \gamma \Vert Im \ z \Vert}, 
\end{equation}
where $\gamma = \min (d, \delta)$,
$a = \max (c_{\delta_1, k}, m_{\delta, k})$ and $b = \max(p, N_1)$. 

Take a number $\gamma_1 \in (0, \gamma)$. 
We can find a convex compact 
$\Omega_k \subset int \ (ch (W_k \cup \overline {\tilde X}_{2, m+1}))$ such that
$$
H_{ch (W_k \cup \overline {\tilde X}_{2, m+1})} (y) - 
H_{\Omega_k}(y) \le \gamma_1 \Vert y \Vert, \ y \in {\mathbb R}^n.
$$
Then from (7) we have that
$$
\vert \widehat {(A^*S_k)}(z) \vert \le a (1 + \Vert z \vert)^b
e^{
H_{\Omega_k + supp \ \mu} (Im \ z)}, 
$$
Note that by the Paley-Wiener-Schwartz theorem \cite [Theorem 7.3.1] {H2} this means that
the support of $A^*S_k$ is contained in $\Omega_k + supp \ \mu$.

Now taking into account the equality 
$$
\widehat {(A^*S_k)}(z) = \hat S_k(z) \hat \mu(z), \ z \in {\mathbb C}^n,
$$
and Proposition 2 we get that $H_{supp \ (A^*S_k)}(x) - H_{supp \ \mu}$ is a supporting function of some convex compact $G_k \subset {\mathbb R}^n$
and for every $\varepsilon > 0$ there exists a constant $C_{\varepsilon} > 0$ such that 
\begin{equation}
\vert \hat S_k(z) \vert \le C_{\varepsilon} \exp (H_{G_k}(Im z) + \varepsilon \Vert z \Vert), \ z \in {\mathbb C}^n.
\end{equation}
Also by the Paley-Wiener-Schwartz theorem \cite [Theorem 7.3.1] {H2} for some $M_k >0$ we have that
$$
\vert \hat S_k(x) \vert \le M_k (1 + \Vert x \Vert)^{N_{2, k}}, \ x \in {\mathbb R}^n,
$$
From this inequality and inequality (8) using Proposition \ref {FL} we get that 
$$
\vert \hat S_k(z) \vert \le M_k (1 + \Vert z \Vert)^{2 N_{2, k}} 
e^{H_{G_k}(Im z)}, \ z \in {\mathbb C}^n.
$$
Again using the Paley-Wiener-Schwartz theorem \cite [Theorem 7.3.1] {H2} we obtain that the support of $S_k$ is contained in $G_k$.
Therefore, for all $y \in {\mathbb R}^n$ we have that
$$
H_{W_k}(y) \le H_{G_k}(y) = H_{supp \ (A^*S_k)}(y) - 
H_{supp \ \mu} (y) \le 
$$
$$
\le H_{\Omega_k + supp \ \mu} (y) - 
H_{supp \ \mu} (y) = H_{\Omega_k} (y).
$$
From this taking into account that $\Omega_k \subset int \ (ch (W_k \cup \overline {\tilde X}_{2, m+1}))$ we get that
$$
H_{W_k}(y) < \max (H_{W_k}(y),
H_{\overline {\tilde X}_{2, m+1}}(y)), \
y \in {\mathbb R}^n.
$$
But it is impossible in view of (\ref {stochka}). Thus, for each $k \in {\mathbb N}$ convex envelope $W_k$ of support of functional $S_k$ $(k=1, 2, \ldots )$ is contained in $\overline {\tilde X}_{2, m+2}$.

Now let $\eta \in {\cal E}({\mathbb R}^n)$ be a function 
with a support in $\overline {\tilde X}_{2, m+4}$ such that 
$0 \le \eta(x) \le 1$ for all $x \in {\mathbb R}^n$ and $\eta(x) = 1$ 
for $x \in \overline {\tilde X}_{2, m+3}$.  
For each $k \in {\mathbb N}$ define a functional $\tilde S_k$ on ${\cal E}(\tilde {X_2})$ by the rule: $\tilde S_k (f) = S_k (\eta f), \ f \in {\cal E}(\tilde {X_2})$. Obviously, $\tilde S_k \in {\cal E}^*(\tilde {X_2})$ and 
$\tilde S_k (f) = S_k (f), \ f \in {\cal E}(X_2)$.
Note that since for each $f \in {\cal E}(X_1)$
$(A+B)(f) = (\tilde {A} + \tilde {B})(f)$
then functionals
$(A+B)^*S_k$ and $(\tilde {A} + \tilde {B})^*{\tilde S_k}$ $(k=1, 2, \ldots )$ coincide
on ${\cal E}(X_1)$. 
Now taking into account that ${\cal E}(X_1)$ is dense in ${\cal E}(\tilde X_1)$  
we get that $(\tilde {A} + \tilde {B})^*{\tilde S_k}$ is the (unique) extension of 
the functional $(A+B)^*S_k$ to ${\cal E}(\tilde X_1)$. 

Show that functionals $(\tilde {A} + \tilde {B})^*{\tilde S_k}$ 
converge in ${\cal E}^*(\tilde X_1)$. 
First note that the sequence $((\tilde {A} + \tilde {B})^*{\tilde S_k})_{k=1}^{\infty}$ is fundamental in ${\cal E}^*(\tilde X_1)$. 
Indeed, let ${\cal B}$ be an arbitrary bounded set in ${\cal E}(\tilde X_1)$ and 
$$
{\cal B}^{\circ} = \{F \in {\cal E}^*(\tilde X_1): 
\vert F(f) \vert \le 1 \ \forall f \in {\cal B} \}
$$
its polar set. Take a function $\omega \in {\cal E}({\mathbb R}^n)$ with a support in $\overline {\tilde X}_{2, m+4} + supp \ \mu$ such that  $0 \le \omega(x) \le 1$ for all 
$x \in {\mathbb R}^n$ and $\omega(x) = 1$ for $x \in \overline {\tilde X}_{2, m+3} + supp \ \mu$. 
Since the support of functionals $\tilde S_k$ is in $\overline {\tilde X}_{2, m+2}$ then 
the support of functionals $(\tilde {A} + \tilde {B})^*{\tilde S_k}$ 
is contained in ${\tilde X}_{2, m+2} + supp \ \mu$. So
for each $f \in {\cal E}(\tilde X_1)$ and all $k, m \in {\mathbb N}$ we have that 
$$
((\tilde {A} + \tilde {B})^*{\tilde S_k}) (f) - ((\tilde {A} + \tilde {B})^*{\tilde S_m}) (f) =((\tilde {A} + \tilde {B})^*{\tilde S_k}) (\omega f) - 
((\tilde {A} + \tilde {B})^*{\tilde S_m}) (\omega f).
$$
We may consider $\omega f$ as an element of ${\cal E}(X_1)$ setting 
$(\omega f) (x) = 0$ for $x \in X_1 \setminus (\overline {\tilde X}_{2, m+4} + supp \ \mu)$.
Then 
$$
((\tilde {A} + \tilde {B})^*{\tilde S_k}) (f) - ((\tilde {A} + \tilde {B})^*{\tilde S_m}) (f) =((A+B)^*S_k) (\omega f) - 
((A+B)^*S_m)(\omega f).
$$
Note that the set $\omega {\cal B} = \{\omega f: f \in {\cal B} \}$ is bounded in 
${\cal E}(X_1)$. Since the sequence $((A+B)^*S_k)_{k=1}^{\infty}$ is converging in 
${\cal E}^*(X_1)$ then it is fundamental in ${\cal E}^*(X_1)$. So there is $N \in {\mathbb N}$ such that for all natural numbers $k, m$: $k, m \ge N$ and $g \in \omega {\cal B}$ we have 
$\vert ((A+B)^*S_k)(g)  - ((A+B)^*S_m)(g) \vert \le 1$. 
Hence, for all natural $k, m$: $k, m \ge N$ and $f \in {\cal B}$ we get that
$$
\vert ((\tilde {A} + \tilde {B})^*{\tilde S_k}) (f) - 
((\tilde {A} + \tilde {B})^*{\tilde S_m}) (f)\vert \le 1.
$$
This means that for all natural $k, m$: $k, m \ge N$ and $f \in {\cal B}$ we have that
$
(\tilde {A} + \tilde {B})^*{\tilde S_k} - 
((\tilde {A} + \tilde {B})^*{\tilde S_m} \in {\cal B}^{\circ}.
$
Thus, we have proved that the sequence $((\tilde {A} + \tilde {B})^*{\tilde S_k})_{k=1}^{\infty}$ is fundamental in ${\cal E}^*(\tilde X_1)$. 
Finally, since ${\cal E}^*(\tilde X_1)$ is complete then we get that 
the sequence $((\tilde {A} + \tilde {B})^*{\tilde S_k})_{k=1}^{\infty}$ 
is converging in ${\cal E}^*(\tilde X_1)$ 
to some element $\tilde T \in {\cal E}^*(\tilde X_1)$. 
But
$(\tilde {A} + \tilde {B})^*({\cal E}^*(\tilde X_2))$
is closed in ${\cal E}^*(\tilde X_1)$.
Hence, there exists a functional $\tilde S \in {\cal E}^*(\tilde X_2)$ such that
$\tilde T = (\tilde {A} + \tilde {B})^*{\tilde S}$.
Let $S$ be restriction of $\tilde S$ on ${\cal E}(X_2)$. Then for each 
$f \in {\cal E}(X_1)$ we have that $\tilde T(f) = T(f)$. Indeed, 
$$
\tilde T(f) = \lim \limits_{k \to \infty}((\tilde {A} + \tilde {B})^*(\tilde S_k))(f)= 
\lim \limits_{k \to \infty}
\tilde S_k((\tilde {A} + \tilde {B}) f) = \lim \limits_{k \to \infty} \tilde S_k ((A + B)f) = 
$$
$$
= \lim \limits_{k \to \infty}  S_k ((A + B)f) = \lim \limits_{k \to \infty} 
((A + B)^*S_k)(f) = T(f).
$$
From this and the following chain of equalities 
$$
\tilde T(f) = \lim \limits_{k \to \infty}((\tilde {A} + \tilde {B})^*(\tilde S_k))(f) 
= ((\tilde {A} + \tilde {B})^*(\tilde S))(f) =
\tilde S((\tilde {A} + \tilde {B}) f) = 
$$
$$
=
\tilde S ((A + B)f) = S((A + B)f) = ((A + B)^*S)(f)
$$
it follows that $T = (A + B)^*S$. 
Thus, the image of the operator $(A+B)^*$ is closed in ${\cal E}^*(X_1)$.
Consequently, the image of the operator $A+B$ is closed in ${\cal E}(X_2)$.

Now we prove that the image of the operator $A+B$ is dense in ${\cal E}(X_2)$. 
It will be done if we show that an arbitrary functional
$S \in {\cal E}^*(X_2)$ with the property
$S((A+B)f) = 0$ for all $f \in {\cal E}(X_1)$ is a zero functional. 
Assume the contrary.
Then the support of $S$ is not empty.
Let $N$ be the order of distribution $S$ and $\delta > 0$ be so small that $(supp \  S)^{\delta} \Subset X_2$. Then there exists a constant $c_{\delta} > 0$  such that  
$$
\vert S(g)\vert \le
c_{\delta}{\Vert g \Vert}_{(supp \  S)^{\delta}, N}, \  g \in {\cal E}(X_2).
$$
From this and the inequality (\ref {ss1}) it follows that there exists a convex compact 
$V \subset int (supp \ \mu)$, a number $N_1 \in {\mathbb Z}_+$ 
(depending on $ch (supp \ S$)) and a constant $C_{\delta} > 0$ 
such that for each $f \in {\cal E}(X_1)$
$$
\vert (B^*S) (f) \vert \le
C_{\delta}{\Vert f \Vert}_{ch (supp \ S) + V, N_1}.
$$
Hence, the support 
of functional $B^*S$ is contained in $ch (supp \   S) + V$.
From the other hand from the equality $B^*S = -\mu*S$ and by the theorem on supports 
\cite[Theorem 4.3.3]{H2}
we have that  $ch (supp \  B^*S) = ch (supp \  S) + supp \  \mu$. 
Thus, $ch (supp \  S) + supp \   \mu \subset ch (supp \   S) + V$. 
But this inclusion is impossible since convex compact $V$ is contained in the interior of the support of $\mu$. 
Hence, our assumption that $S$ is not a zero functional was false. 
Thus, $S=0$. This means that the image of the operator $A+B$ is dense in ${\cal E}(X_2)$.

Theorem is proved.

\begin{section}
{\bf  Example of the operator $B$} 
\end{section}

Let $\mu \in {\cal E}'({\mathbb R}^n)$ be an invertible distribution and
$supp \ \mu = \overline {D(1)}$. Distributions with these properties can be constructed 
(see, e.g., \cite [Theorem 1, Theorem 3, Theorem 4] {A}).
Let $X_2=D(1)$, $X_1 = D(2)$.
Let
$A: {\cal E}(X_1	) \to {\cal E}(X_2)$ be a convolution operator acting by the rule
$
(Af)(x) = (\mu * f)(x), \ x \in X_1.
$
Take a function $b \in {\cal E}({\mathbb R}^{2n})$ with the support in 
$\overline {D(\frac {1}{4})} \times \overline {D(\frac {1}{4})}$.
Define the operator
$B: {\cal E}(X_1) \to {\cal E}(X_2)$ acting by the rule
$$
(Bf)(x) =  \int_{{\mathbb R}^n} b(x, \xi) f(x + \xi) \ d \xi, \ \Vert x \Vert \le \frac {1}{4},
$$
$$
(Bf)(x) = 0, \ \frac {1}{4} < \Vert x \Vert < 1.
$$

Let $K$ be a convex compact in $X_2$ and $\gamma: = dist (K,  \partial X_2)$. 
Show that there exists a convex compact 
$V \subset int (supp \ \mu)$ such that for any 
$\varepsilon \in (0, \gamma)$ and for any $N_2 \in {\mathbb Z}_+$ there exists a constant 
$c = c(\varepsilon, N_2) > 0 $ such that
$$
\label{s1}
{\Vert Bf \Vert}_{K^{\varepsilon}, N_2} \le
c {\Vert f \Vert}_{K + V,  0} \  , \
f \in {\cal E}(X_1).
$$
Obviously, for each $\varepsilon \in (0, \gamma)$
and for each $N_2 \in {\mathbb Z}_+$ there exists a constant $C > 0$ depending on 
$b$ and $N_2$ such that for each $f \in {\cal E}(X_1)$ we have that
\begin {equation}
{\Vert Bf \Vert}_{K^{\varepsilon}, N_2} = 
{\Vert Bf \Vert}_{K^{\varepsilon} \cap \overline {D(\frac {1}{4})}, N_2} \le 
C_1 \Vert f \Vert_{(K^{\varepsilon} \cap \overline {D(\frac {1}{4})}) + 
\overline {D(\frac {1}{4})} ,  0}.
\end {equation}
If $\gamma \in (0, \frac {3}{4})$ then from (9) we have that
$$
{\Vert Bf \Vert}_{K^{\varepsilon}, N_2} \le C_1 \Vert f \Vert_{K^{\gamma} + 
\overline {D(\frac {1}{4})} , 0} =  C_1 \Vert f \Vert_{K + 
\overline {D(\gamma + \frac {1}{4})} , 0}
$$
Hence, in this case we can put $V = \overline {D(\gamma + \frac {1}{4})}$. 
If $\gamma \in [\frac {3}{4}, 1]$ then $K \subset \overline {D(\frac {1}{4})}$ and from (9) we have that 
$$
{\Vert Bf \Vert}_{K^{\varepsilon}, N_2} 
\le
C_1 \Vert f \Vert_{\overline {D(\frac {1}{2})},  0} \le 
C_1 \Vert f \Vert_{K + 
\overline {D(\frac {3}{4})} , 0}.
$$
So if $\gamma \in [\frac {3}{4}, 1]$ then we can put $V = \overline {D(\frac {3}{4})}$. 

Thus, by Theorem the operator $A+B: {\cal E}(X_1) \to {\cal E}(X_2)$ is surjective. 

{\bf Acknowledgements}. The research was supported by grant from RFBR (15-01-01661) and RAS Program for Fundamental Research (Modern Problems of Theoretical Mathematics, project "Complex Analysis and Functional Equations").

\end{document}